\documentclass[11pt]{article}
\usepackage{amsmath}
\usepackage{amssymb}
\usepackage{amscd}
\usepackage{verbatim}
\newcommand{\be}{\begin{equation}}
\newcommand{\bs}{\begin{sub}}
\newcommand{\es}{\end{sub}}
\newcommand{\bsn}{\begin{subn}}
\newcommand{\esn}{\end{subn}}
\newcommand{\bea}{\begin{eqnarray}}
\newcommand{\eea}{\end{eqnarray}}
\newcommand{\bean}{\begin{eqnarray*}}
\newcommand{\eean}{\end{eqnarray*}}
\newcommand{\BA}[1]{\begin{array}{#1}}
\newcommand{\EA}{\end{array}}
\newcommand{\Real}{\mathbb{R}}

\newlength{\wex}  \newlength{\hex}
            \def\gl{\lambda}

\newcommand{\pf}{\noindent \mbox{{\bf Proof}: }}
\def\squarebox#1{\hbox to #1{\hfill\vbox to #1{\vfill}}}
\newcommand{\qed}{\hspace*{\fill}
      \vbox{\hrule\hbox{\vrule\squarebox{.667em}\vrule}\hrule}\smallskip}

\newcommand{\Green}[4]{\mbox{$G^{#1}_{#2}(#3,#4)$}}

\newcommand{\beq}{\begin{equation}}
\newcommand{\eeq}{\end{equation}}
\newcommand{\beqa}{\begin{eqnarray}}
\newcommand{\eeqa}{\end{eqnarray}}
\newcommand{\beqanl}{\begin{eqnarray*}}
\newcommand{\eeqanl}{\end{eqnarray*}}
\begin{document}
\renewcommand{\theequation}{\thesection.\arabic{equation}}
\newcommand{\mysection}[1]{\section{#1}\setcounter{equation}{0}}

\def\stackunder#1#2{\mathrel{\mathop{#2}\limits_{#1}}}
\newtheorem{theorem}{Theorem}[section]
\newtheorem{lemma}[theorem]{Lemma}
\newtheorem{definition}[theorem]{Definition}
\newtheorem{corollary}[theorem]{Corollary}
\newtheorem{conjecture}[theorem]{Conjecture}
\newtheorem{remark}[theorem]{Remark}
\newtheorem{question}[theorem]{Question}
\newtheorem{example}[theorem]{Example}
\newtheorem{Thm}[theorem]{Theorem}
\newtheorem{Lem}[theorem]{Lemma}
\newtheorem{Pro}[theorem]{Proposition}
\newtheorem{Def}[theorem]{Definition}
\newtheorem{Exa}[theorem]{Example}
\newtheorem{Exs}[theorem]{Examples}
\newtheorem{Rems}[theorem]{Remarks}
\newtheorem{Rem}[theorem]{Remark}
\newtheorem{Cor}[theorem]{Corollary}
\newtheorem{Conj}[theorem]{Conjecture}
\newtheorem{Prob}[theorem]{Problem}
\newtheorem{Ques}[theorem]{Question}

\def\stackunder#1#2{\mathrel{\mathop{#2}\limits_{#1}}}

\author{Yehuda Pinchover
\\ Department of Mathematics\\  Technion - Israel Institute of Technology\\
Haifa 32000, Israel\\ pincho@techunix.technion.ac.il}
\title{Large time behavior of the heat kernel}
\date{}
\maketitle
\begin{abstract} In this paper we study the large time behavior of the (minimal)
heat kernel $k_P^M(x,y,t)$ of a general time independent parabolic
operator $L=u_t+P(x, \partial_x)$ which is defined on a noncompact
manifold $M$. More precisely, we prove that $$\lim_{t\to\infty}
e^{\lambda_0 t}k_P^{M}(x,y,t)$$ always exists. Here $\lambda_0$ is
the generalized principal eigenvalue of the operator $P$ in $M$.
\\[1mm]

\end{abstract}
\mysection{Introduction}\label{Introduction}
    Let $k_P^{M}(x,y,t)$ be the (minimal) heat kernel of a time independent
parabolic operator $Lu=u_t+P(x, \partial_x)u$ which is defined on
a noncompact Riemannian manifold $M$. Denote by $\lambda_0$ the
generalized (Dirichlet) principal eigenvalue of the operator $P$
in $M$.

Over the past three decades, there have been a large number of
works devoted to large time estimates of the heat kernel in
various settings (see for example the following monographs and
survey articles
\cite{Ch,C,Dheat,G,Has,Pinsky,PE,R,Ssemi,Va,VSCC,W}, and the
references therein). Despite the wide diversity of the results in
this field, the following basic question has not been fully
answered.
\begin{question}\label{qu}
Does  $\lim_{t\to\infty} e^{\lambda_0 t}k_P^{M}(x,y,t)$ always
exist?
\end{question}
The aim of this paper is to give a complete answer to Question
\ref{qu} for arbitrary $P$ and $M$.  The following theorem
\cite{Pheat} gives only a partial answer to the above question
(see also \cite{CK,Has,Pinsky,S,Va}).
\begin{theorem}\label{thm1} Let $P$ an elliptic operator on $M$.
\begin{description}
\item[{\em (i)}] If $P-\lambda_0$ is subcritical in $M$ (i.e. $\int_0^\infty
e^{\lambda_0 t}k_P^{M}(x,y,t)dt<\infty$), then $$\lim_{t\to\infty}
e^{\lambda_0 t} k_P^{M}(x,y,t)=0.$$
  \item[{\em (ii)}] If $P-\lambda_0$ is positive-critical in $M$ (i.e.
$\int_0^\infty e^{\lambda_0 t}k_P^{M}(x,y,t)\,dt=\infty$, and the
ground states $\varphi$ and $\varphi^*$ of $P-\lambda_0$ and
$P^*-\lambda_0$ respectively, satisfy $\varphi^*\varphi \in
L^1(M)$), then
  $$\lim_{t\to\infty} e^{\lambda_0 t} k_P^{M}(x,y,t)=
  \frac{\varphi(x)\varphi^*(y)}{\int_M \varphi^*(z)\varphi(z)\,dz}\,.$$
\item[{\em (iii)}] If $P-\lambda_0$ is null-critical in $M$ (i.e.
$\int_0^\infty e^{\lambda_0 t}k_P^{M}(x,y,t)\,dt=\infty$, and the
ground states $\varphi$ and $\varphi^*$ of $P-\lambda_0$ and
$P^*-\lambda_0$ respectively, satisfy $\varphi^*\varphi\not \in
L^1(M)$), then $$\lim_{T\to\infty}
  \frac{1}{T}\int_0^T e^{\lambda_0 t}k_P^{M}(x,y,t)\,dt=0.$$
Moreover, if one assumes further that $P$ is a formally symmetric
operator ($P=P^*$), then in the null-critical case
$\lim_{t\to\infty}
  e^{\lambda_0 t} k_P^{M}(x,y,t)=0$.
\end{description}
\end{theorem}

The main result of the present paper is the following theorem
which answers the author's conjecture \cite[Remark 1.4]{Pheat}
about the existence of the limit in the null-critical nonsymmetric
case.
\begin{theorem}\label{mainthm}
Assume that $P-\lambda_0$  is a (nonsymmetric) null-critical
operator in $M$.  Then $$\lim_{t\to\infty} e^{\lambda_0 t}
k_P^{M}(x,y,t)=0.$$
\end{theorem}
Thus, theorems  \ref{thm1} and \ref{mainthm} indeed solve Question
\ref{qu}. More precisely, with the aid of Theorem 1.1 of
\cite{Pheat}, we have
\begin{corollary}\label{cormain}
The $\lim_{t\to\infty} e^{\lambda_0 t}k_P^{M}(x,y,t)$ exists for
all $x,y\in M$, and the limit is positive if and only if the
operator $P-\lambda_0$ is positive-critical.

Moreover, let $\Green{M}{P-\lambda}{x}{y}$ be the minimal positive
Green function of the elliptic operator $P-\lambda$ on $M$. Then
\begin{equation}\label{eqgreen}
\lim_{t\to\infty} e^{\lambda_0 t}k_P^{M}(x,y,t)=
\lim_{\lambda\nearrow\lambda_0}(\lambda_0-\lambda)\Green{M}{P-\lambda}{x}{y}.
\end{equation}
\end{corollary}

\vspace{3mm}

The proof of Theorem \ref{mainthm} hinges on Lemma \ref{lemVa}
which is a slight extension of a lemma of Varadhan (see,
\cite[Lemma 9, page 259]{Va} or \cite[pp. 192--193]{Pinsky}).
Varadhan proved his lemma for positive-critical operators on
$\mathbb{R}^d$ using a purely probabilistic approach. Our key
observation is that the assertion of Varadhan's lemma is valid
under the weaker assumption that the skew product operator
$\bar{P}= P\otimes I+I\otimes P$ is critical in $\bar{M}=M\times
M$, where $I$ is the identity operator on $M$. We note that if
$\bar{P}$ is subcritical in $\bar{M}$, then by Theorem \ref{thm1},
the heat kernel of $\bar{P}$ on $\bar{M}$ tends to zero as
$t\to\infty$. Since the heat kernel of $\bar{P}$ is equal to the
product of the heat kernels of its factors, it follows that if
$\bar{P}$ is subcritical in $\bar{M}$, then $\lim_{t\to\infty}
k_P^{M}(x,y,t)=0$.

In Section \ref{secVa}, we formulate and give a purely analytic
proof of Lemma \ref{lemVa}. Our proof of the lemma is in fact the
translation of Varadhan's proof to the analytic apparatus. It uses
the large time behaviors of the parabolic capacitory potential and
of the heat content (see Section \ref{secaux}).

The proof of Theorem \ref{mainthm} is given in Section
\ref{secmainthm}. We conclude the paper with some open problems
which are closely related to the large time behavior of the heat
kernel (see Section \ref{secopen}).
\begin{remark}\label{remdec}{\em
In the null-recurrent case, the heat kernel may decay very slowly
as $t\to \infty$, and one can construct a complete Riemannian
manifold $M$ such that all its Riemannian products $M^j, j\geq 1$
are null-recurrent (see \cite{CG}).
 }\end{remark}
\begin{remark}
{\em We would like to point out that the results of this paper,
are also valid for an elliptic operator $P$ in divergence form and
also for a strongly elliptic operator $P$ with locally bounded
coefficients.}
\end{remark}
\mysection{Preliminaries}\label{sec1}
 Let $P$ be a linear, second
order, elliptic operator defined in a  noncompact, connected,
$C^3$-smooth Riemannian manifold $M$ of dimension $d$. Here $P$ is
an  elliptic operator with real, H\"{o}lder continuous
coefficients which in any coordinate system
$(U;x_{1},\ldots,x_{d})$ has the form \be \label{P}
P(x,\partial_{x})=-\sum_{i,j=1}^{d}
a_{ij}(x)\partial_{i}\partial_{j} + \sum_{i=1}^{d}
b_{i}(x)\partial_{i}+c(x),
\end{equation}
where $\partial_{i}=\partial/\partial x_{i}$. We assume that for
every $x\in M$ the real quadratic form
\be
\label{ellip} \sum_{i,j=1}^{d} a_{ij}(x)\xi_{i}\xi_{j},\;\;
\xi=(\xi_{1},\ldots,\xi_{d}) \in \Real ^d
\end{equation}
is positive definite. The formal adjoint of $P$ is denoted by
$P^*$. We consider the parabolic operator $L$
\begin{equation}\label{eqL}
  Lu=u_t+Pu \qquad \mbox{ on } M\times (0,\infty).
\end{equation}

Let $\{M_{j}\}_{j=1}^{\infty}$ be an {\em exhaustion} of  $M$,
i.e. a sequence of smooth, relatively compact domains such that
$M_1\neq \emptyset$, $\mbox{cl}({M}_{j})\subset M_{j+1}$ and
$\cup_{j=1}^{\infty}M_{j}=M$. For every $j\geq 1$, we denote
$M_{j}^*=M\setminus \mbox{cl}({M_j})$. Let
$M_\infty=M\cup\{\infty\}$ be the one-point compactification of
$M$. By the notation $x\to\infty$, we mean that $x\to\infty$ in
the topology of $M_\infty$.

Denote the cone of all positive (classical) solutions of the
equation $Pu\!=\!0$ in $M$ by $\mathcal{C}_{P}(M)$. The {\em
generalized principal eigenvalue}  is defined by
$$\gl_0=\gl_0(P,M)
:= \sup\{\gl \in \mathbb{R} \; :\; \mathcal{C}_{P-\lambda}(M)\neq \emptyset\}.$$
Throughout this paper we always assume that $\lambda_0\geq 0$.

For every $j\!\geq\! 1$, consider the Dirichlet heat kernel
$k_P^{M_j}(x,y,t)$ of the parabolic operator $L=\partial_t+P$ in
$M_j$. So, for every continuous function $f$ with a compact
support in $M$, $u(x,t)=\int_{M_j} k_P^{M_j}(x,y,t)f(y)\, dy$
solves the initial-Dirichlet boundary value problems
 \bea\label{eqibvpj}
 Lu&=&0 \quad \mbox{ in } M_j\times (0,\infty),\nonumber\\
  u&=&0 \quad \mbox{ on  } \partial M_j\times (0,\infty),\\
  u&=&f  \quad \mbox{ on  }  M_j\times \{0\}.\nonumber
  \eea
   By the maximum principle, $\{k_P^{M_j}(x,y,t)\}_{j=1}^{\infty}$ is
an increasing sequence which converges to $k_P^{M}(x,y,t)$, the
{\em minimal heat kernel} of the parabolic operator $L$ in $M$. If
$$\int_0^\infty k_P^{M}(x,y,t)\,dt<\infty \qquad
\mbox{(respectively, $\int_0^\infty
k_P^{M}(x,y,t)\,dt=\infty$),}$$ then $P$ is said to be a {\em
subcritical} (respectively, {\em critical}) operator in $M$,
\cite{Pinsky}.

It can be easily checked that for $\lambda\leq \lambda_0$, the
heat kernel $k_{P-\lambda}^M$ of the operator $P-\lambda$ is equal
to $e^{\lambda t}k_P^M(x,y,t)$. Since we are interested in the
asymptotic behavior of $e^{\lambda_0 t}k_P^M(x,y,t)$, we assume
throughout the paper (unless otherwise stated) that $\lambda_0=0$.

It is well known that if $\lambda_0>0$, then $P$ is subcritical in
$M$. Clearly, $P$ is critical (respectively, subcritical) in $M$,
if and only if $P^*$ is critical (respectively, subcritical) in
$M$. Furthermore, if $P$ is critical in $M$, then
$\mathcal{C}_{P}(M)$ is a one-dimensional cone. In this case,
$\varphi \in \mathcal{C}_{P}(M)$ is called a {\em ground state of
the operator $P$ in $M$} \cite{Pheat,Pinsky}. We denote the ground
state of $P^*$ by $\varphi^*$.

The ground state $\varphi$ is a global positive solution of the
equation $Pu=0$ of {\em minimal growth in a neighborhood of
infinity in} $M$. That is, if $v\in C(M_j^*)$ is a positive
solution of the equation $Pu=0$ in $M_j^*$ such that $\varphi\leq
v$ on $\partial M_j^*$, then $\varphi\leq v$ in $M_j^*$
\cite{Pheat,Pinsky}.

In the critical case, the ground state $\varphi$ (respectively,
$\varphi^*$) is a positive invariant solution of the operator $P$
(respectively, $P^*$) in $M$ (see for example
\cite{Pheat,Pinsky}). That is,
\begin{equation}\label{eqinvar}
\int_M k_P^M(x,y,t)\varphi(y)\, dy= \varphi(x),\;\;  \mbox{and }
\int_M k_P^M(x,y,t)\varphi^*(x)\, dx= \varphi^*(y).
\end{equation}

\begin{definition}\label{defnull}{\em
A critical operator $P$ is said to be {\em positive-critical} in
$M$ if $\varphi^*\varphi\in L^1(M)$, and {\em null-critical} in
$M$ if $\varphi^*\varphi\not\in L^1(M)$.
 }\end{definition}
\begin{remark}{\em Let $\mathbf{1}$ be the constant function on $M$,
taking at any point $x\in M$ the value $1$. Suppose that
$P\mathbf{1}=0$. Then $P$ is subcritical (respectively,
positive-critical, null-critical) in $M$ if and only if the
corresponding diffusion process is transient (respectively,
positive-recurrent, null-recurrent) \cite{Pinsky}.
 In fact, since we are interested in the critical case, it is natural to use the
$h$-transform with $h=\varphi$. So, $$P^\varphi
u=\frac{1}{\varphi}P(\varphi u)\qquad\mbox{ and }\quad
k_{P^\varphi}^M(x,y,t)=\frac{1}{\varphi(x)}k_{P}^M(x,y,t)\varphi(y).$$
Note that $P^\varphi$ is null-critical (respectively,
positive-critical) if and only if $P$ is null-critical
(respectively, positive-critical), and the ground states of
$P^\varphi$ and $(P^\varphi)^*$ are $\mathbf{1}$ and
$\varphi^*\varphi$, respectively. Moreover,
$$\lim_{t\to\infty}k_{P^\varphi}^M(x,y,t)=0 \quad \mbox{ if and
only if } \quad
 \lim_{t\to\infty}k_{P}^M(x,y,t)=0.$$ Therefore, throughout the
paper (unless otherwise stated), we assume that $$\mbox{{\bf
(A)}}\qquad \qquad\qquad P\mathbf{1}=0, \mbox{ and } P \mbox{ is a
critical operator in } M.\qquad\qquad\qquad\mbox{}$$
 }\end{remark}

It is well known that on a general noncompact manifold $M$, the
solution of the Cauchy problem for the parabolic equation $Lu=0$
is not uniquely determined (see for example \cite{IM} and the
references therein). On the other hand, under Assumption {\bf
(A)}, there is a unique {\em minimal} solution of the Cauchy
problem and of certain initial-boundary value problems for {\em
bounded} initial and boundary conditions. More precisely,
\begin{definition}\label{defCp}{\em
Let $f$ be a bounded continuous function on $M$.  By the {\em
minimal solution} $u$ of the Cauchy problem \bean\label{eqCp}
 Lu&=&0 \quad \mbox{ in } M\times (0,\infty),\\
  u&=&f  \quad \mbox{ on  }  B^*\times \{0\},
  \eean
we mean the function
\begin{equation}\label{eqminim}
u(x,t):=\int_{M} k_P^{M}(x,y,t)f(y)\,dy.
 \end{equation}
 }\end{definition}
\begin{definition}\label{defibvp}{\em
Let $B\subset\subset M$ be a smooth bounded domain such that
$B^*:=M\setminus \mbox{cl}(B)$ is connected. Assume that $f$ is a
bounded continuous function on $B^*$, and $g$ is a bounded
continuous function on $\partial B\times (0,\infty)$.  By the {\em
minimal solution} $u$ of the initial-boundary value problem
 \bea\label{eqibvpfg}
 Lu&=&0 \quad \mbox{ in } B^*\times (0,\infty),\nonumber\\
  u&=&g \quad \mbox{ on  } \partial B\times (0,\infty),\\
  u&=&f  \quad \mbox{ on  }  B^*\times \{0\}\nonumber,
  \eea
we mean the limit of the solutions $u_j$ of the following
initial-boundary value problems \bean\label{eqibvpfgk}
 Lu&=&0 \quad \mbox{ in } (B^*\cap M_j)\times (0,\infty),\\
  u&=&g \quad \mbox{ on  } \partial B\times (0,\infty),\\
    u&=&0 \quad \mbox{ on  } \partial M_j\times (0,\infty),\\
  u&=&f  \quad \mbox{ on  }  (B^*\cap M_j)\times \{0\}.
  \eean
}\end{definition}
\begin{remark}{\em
It can be easily checked that the sequence $\{u_j\}$ is indeed a
converging sequence which converges  to a solution of the
initial-boundary value problem (\ref{eqibvpfg}).
 }\end{remark}
 \mysection{Auxiliary results}\label{secaux}
\label{Auxiliary results}
\begin{lemma}\label{lem3}
 Assume that $P\mathbf{1}=0$ and that $P$ is critical in $M$. Let
$B:=B(x_0,\delta)\subset\subset M$ be the ball of radius $\delta$
centered at $x_0$, and suppose that $B^*=M\setminus \mbox{cl}(B)$
is connected. Let $w$ be {\bf the heat content of $B^*$}, i.e. the
minimal nonnegative solution of the following initial-boundary
value problem
 \bea\label{eqw3}
 Lu&=&0 \quad \mbox{ in } B^*\times (0,\infty),\nonumber\\
  u&=&0 \quad \mbox{ on  } \partial B\times (0,\infty),\\
  u&=&1  \quad \mbox{ on  }  B^*\times \{0\}.\nonumber
  \eea
Then $w$ is a decreasing function of $t$, and
$\lim_{t\to\infty}w(x,t)=0$ locally uniformly in $M$.
\end{lemma}
\pf Clearly,
\begin{equation}\label{eqw4}
w(x,t)=\int_{B^*} k_P^{B^*}(x,y,t)\,dy<\int_{M}
k_P^{M}(x,y,t)\,dy=1 .
\end{equation}
 It follows that $0<w<1$ in $B^*\times (0,\infty)$. Let
 $\varepsilon>0$.
By the semigroup identity and (\ref{eqw4}),
 \bea\label{eqw5}
w(x,t+\varepsilon)=\int_{B^*}
k_P^{B^*}(x,y,t+\varepsilon)\,dy&=&\nonumber\\[2mm] \int_{B^*}
\left(\int_{B^*}k_P^{B^*}(x,z,t)k_P^{B^*}(z,y,\varepsilon)\,dz\right)
\,dy&=&\\[2mm]
\int_{B^*}k_P^{B^*}(x,z,t)\left(\int_{B^*}k_P^{B^*}(z,y,\varepsilon)\,dy
\right)\,dz&<&\int_{B^*}k_P^{B^*}(x,z,t)\,dz =w(x,t).\nonumber
\eea Hence, $w$ is a decreasing function of $t$, and therefore,
$\lim_{t\to\infty}w(x,t)$ exists.

We denote the limit function by $v$. So, $0\leq v<\mathbf{1}$ and
$v$ is a solution of the elliptic equation $Pu=0$ in $B^*$ which
satisfies $u=0$  on $\partial B$. Therefore, $\mathbf{1}-v$ is a
positive solution of the equation $Pu=0$ in $B^*$ which satisfies
$u=1$  on $\partial B$. On the other hand, it follows from the
criticality assumption that $\mathbf{1}$ is the minimal positive
solution of the equation $Pu=0$ in $B^*$ which satisfies $u=1$  on
$\partial B$. Thus, $\mathbf{1}\leq \mathbf{1}-v$, and therefore,
$v=0$.
 \qed
\begin{definition}\label{defcapac}
{\em Let $B:=B(x_0,\delta)\subset\subset M$. Suppose that
$B^*=M\setminus \mbox{cl}(B)$ is connected. The nonnegative
(minimal) solution $$v(x,t)=\mathbf{1}-\int_{B^*}
k_P^{B^*}(x,y,t)\,dy$$ is called the {\em parabolic capacitory
potential of $B^*$}. Note that $v$ is indeed the minimal
nonnegative solution of the initial-boundary value problem
 \bea\label{fet}
 Lu&=&0 \quad \mbox{ in } B^*\times (0,\infty),\nonumber\\
  u&=&1 \quad \mbox{ on  } \partial B\times (0,\infty),\\
  u&=&0  \quad \mbox{ on  }  B^*\times \{0\}.\nonumber
  \eea
 }\end{definition}

\begin{corollary}\label{lem2}
Under the assumptions of Lemma \ref{lem3},  the parabolic
capacitory potential $v$ of $B^*$  is an increasing function of
$t$, and $\lim_{t\to\infty}v(x,t)=1$ locally uniformly in $M$.
\end{corollary}
\pf Clearly,
\begin{equation}\label{eqv1}
v(x,t)=\mathbf{1}-\int_{B^*}
k_P^{B^*}(x,y,t)\,dy=\mathbf{1}-w(x,t)
\end{equation}
where $w$ is the heat content of $B^*$. Therefore, the corollary
follows directly from Lemma \ref{lem3}.\qed
\mysection{Varadhan's lemma}\label{secVa} In this section, we give
a purely analytic proof of a lemma of Varadhan \cite[Lemma 9, page
259]{Va} for a slightly more general case. We consider the
Riemannian product manifold $\bar{M}:=M\times M$. A point in
$\bar{M}$ is denoted by $\bar{x}=(x_1,x_2)$. Let $P_{x_i}$,
$i=1,2$ denote the operator $P$ in the variable $x_i$, and let
$\bar{P}=P_{x_1}+P_{x_2}$ be the skew product operator defined on
$\bar{M}$. We denote by $\bar{L}$ the corresponding parabolic
operator. Note that if $\bar{P}$ is critical in $\bar{M}$, then
$P$ is critical in $M$. Moreover, if $P$ is positive-critical in
$M$, then $\bar{P}$ is positive-critical in $\bar{M}$.
 \begin{lemma}\label{lemVa} Assume that $P\mathbf{1}=0$.
Suppose further that $\bar{P}$ is critical on $\bar{M}$. Let $f$
be a continuous bounded function on $M$, and  let
 $$u(x,t)=\int_{M} k_P^{M}(x,y,t)f(y)\,dy$$
 be the minimal solution of the Cauchy problem with initial data $f$ on
$M$.  Fix $K\subset \subset M$. Then
 $$\lim_{t\to\infty}\sup_{x_1,x_2\in K}|u(x_1,t)-u(x_2,t)|=0.$$
\end{lemma}
\pf Denote by $\bar{u}(\bar{x},t):=u(x_1,t)-u(x_2,t)$. Recall that
the heat kernel $\bar{k}(\bar{x},\bar{y},t)$ of the operator
$\bar{L}$ on $\bar{M}$ satisfies
\begin{equation}\label{eqkkk}
\bar{k}_{\bar{P}}^{\bar{M}}(\bar{x},\bar{y},t)=k_P^M(x_1,y_1,t)k_P^M(x_2,y_2,t).
\end{equation}
By (\ref{eqinvar}) and (\ref{eqkkk}), we have
 \bean
\bar{u}(\bar{x},t)=u(x_1,t)-u(x_2,t)=\\
\int_{M}k_P^M(x_1,y_1,t)f(y_1)\,dy_1-\int_{M}k_P^M(x_2,y_2,t)f(y_2)\,dy_2=\\
\int_{M}\int_{M}k_P^M(x_1,y_1,t)k_P^M(x_2,y_2,t)(f(y_1)-f(y_2))\,dy_1dy_2=\\
\int_{\bar{M}}\bar{k}_{\bar{P}}^{\bar{M}}
(\bar{x},\bar{y},t)(f(y_1)-f(y_2))\,d\bar{y}.
 \eean
Hence, $\bar{u}$ is the minimal solution of the Cauchy problem for
the equation $\bar{L}\bar{u}=0$ with initial data $f(x_1)-f(x_2)$
on $\bar{M}$.

Fix a compact set $K\subset\subset M$ and $x_0\in M\setminus K$,
and let $\varepsilon>0$. Let $B:=B((x_0,x_0),\delta)\subset\subset
\bar{M}\setminus \bar{K}$, where $\bar{K}=K\times K$, and $\delta$
will be determined below. We may assume that $B^*=\bar{M}\setminus
\mbox{cl}(B)$ is connected. Then $\bar{u}$ is a minimal solution
of the following initial-boundary value problem
 \bea\label{ibvpf}
 \bar{L}\bar{u}&=&0 \qquad \qquad\qquad\qquad \mbox{ in }  B^*\times (0,\infty),\nonumber\\
  \bar{u}(\bar{x},t)&=&u(x_1,t)-u(x_2,t) \quad
  \mbox{ on } \partial B\times (0,\infty),\\
  \bar{u}(\bar{x},0)&=&f(x_1)-f(x_2)  \qquad\quad  \mbox{ on  }  B^*\times \{0\}.\nonumber
  \eea
We need to prove that $\lim_{t\to\infty}\bar{u}(\bar{x},t)=0$.

By the superposition principle (which obviously holds for minimal
solutions), we have
 $$\bar{u}(\bar{x},t)=u_1(\bar{x},t)+u_2(\bar{x},t)\quad \mbox{ on } B^*\times
[1,\infty),$$ where $u_1$ solves the initial-boundary value
problem
 \bea\label{ibvpu1}
 \bar{L}u_1&=&0 \qquad\qquad\qquad\qquad \mbox{ in } B^*\times (1,\infty),\nonumber\\
  u_1(\bar{x},t)&=&u(x_1,t)-u(x_2,t) \quad \mbox{ on  } \partial B\times (1,\infty),\\
  u_1(\bar{x},0)&=&0  \qquad\qquad\qquad\qquad \mbox{ on  }  B^*\times \{1\},\nonumber
 \eea
 and $u_2$ solves the  initial-boundary value problem
 \bea\label{ibvpu2}
 \bar{L}u_2&=&0 \qquad\qquad\qquad\qquad \mbox{ in } B^*\times (1,\infty),\nonumber\\
  u_2(\bar{x},t)&=&0 \qquad\qquad\qquad\qquad \mbox{ on  } \partial B\times (1,\infty),\\
  u_2(\bar{x},0)&=& u(x_1,1)-u(x_2,1)  \quad \mbox{ on  }  B^*\times \{1\}.\nonumber
  \eea
Clearly, $|\bar{u}(\bar{x},t)|\leq 2\|f\|_\infty$ on
$\bar{M}\times (0,\infty)$. Note that if $\bar{x}=(x_1,x_2)\in
\partial B$, then on $M$, $\mbox{dist}_M(x_1,x_2) < 2\delta$. Using Schauder's
parabolic interior estimates on $M$, it follows that if $\delta$
is small enough, then
$$|\bar{u}(\bar{x},t)|=|u(x_1,t)-u(x_2,t)|<\varepsilon \quad
\mbox{ on  }
\partial B\times (1,\infty).$$
By comparison of $u_1$ with the parabolic capacitory potential of
$B^*$, we obtain that
\begin{equation}\label{estu1}
|u_1(\bar{x},t)|\leq \varepsilon\left(1-\int_{B^*}
\bar{k}_{\bar{P}}^{B^*}(\bar{x},\bar{y},t-1)\,
d\bar{y}\right)<\varepsilon \qquad \mbox{ in } B^*\times
(1,\infty).
  \end{equation}
On the other hand,
\begin{equation}\label{estu2}
|u_2(\bar{x},t)|\leq 2\|f\|_\infty\int_{B^*}
\bar{k}_{\bar{P}}^{B^*}(\bar{x},\bar{y},t-1)\,d\bar{y} \qquad
\mbox{ in } B^*\times (1,\infty).
  \end{equation}
It follows from (\ref{estu2}) and Lemma \ref{lem3} that there
exists $T>0$ such that
\begin{equation}\label{estu21}
|u_2(\bar{x},t)|\leq \varepsilon \quad \mbox { for all }
\bar{x}\in \bar{K} \mbox{ and } t>T.
\end{equation}
Combining (\ref{estu1}) and (\ref{estu21}), we obtain that
$|u(x_1,t)-u(x_2,t)|\leq 2\varepsilon$ for all $x_1,x_2\in K $ and
$t>T$. Since $\varepsilon$ is arbitrary, the lemma is proved. \qed

\mysection{Proof of Theorem \ref{mainthm}}\label{secmainthm}
Without loss of generality, we may assume that $P\mathbf{1}=0$,
where $P$ is a null-critical operator in $M$. We need to prove
that
 $\lim_{t\to\infty} k_P^M(x,y,t)=0$.

Consider again the Riemannian product manifold $\bar{M}:=M\times
M$ and let $\bar{P}=P_{x_1}+P_{x_2}$ be the corresponding skew
product operator which is defined on $\bar{M}$. If $\bar{P}$ is
subcritical on $\bar{M}$, then by Theorem \ref{thm1},  $\lim_{t\to
\infty} \bar{k}_{\bar{P}}^{\bar{M}}(x,y,t)=0$. Since
$$\bar{k}_{\bar{P}}^{\bar{M}}(\bar{x},\bar{y},t)=
k_{P}^M(x_1,y_1,t)k_{P}^M(x_2,y_2,t),$$ it follows that
$\lim_{t\to \infty}k_P^M(x,y,t)=0$.

Therefore, there remains to prove the theorem for the case where
$\bar{P}$ is critical in $\bar{M}$. Fix a nonnegative, bounded,
continuous function $f\neq 0$ such that $\varphi^*f\in L^1(M)$,
and consider the solution
 $$v(x,t)=\int_M k_P^M(x,y,t)f(y)\,dy.$$
Let $t_n\to \infty$. then by subtracting a subsequence, we may
assume that for any $t\in \mathbb{R}$ the function $v(x,t+t_n)$
converges to a nonnegative solution $u\in \mathcal{H}_+(M\times
\mathbb{R})$, where $$\mathcal{H}_+(M\times\mathbb{R}):=\{u\geq
0\, |\,Lu=0 \mbox{ in } M\times \mathbb{R}\}.$$

Invoking Lemma \ref{lemVa} (Varadhan's lemma), we see that
$u(x,t)=\alpha(t)$. Since $u$ solves the parabolic equation
$Lu=0$, it follows that $\alpha(t)$ is a nonnegative constant
$\alpha$.

We claim that $\alpha=0$. Suppose to the contrary that $\alpha>0$.
The assumption that $\varphi^*f\in L^1(M)$ and (\ref{eqinvar})
imply that for any $t>0$
 \bea\label{eqphiv} \int_M
\varphi^*(y)v(y,t)\,dy=\int_M \varphi^*(y)\left(\int_M
k_P^M(y,z,t)f(z)\,dz\right)\,dy=\nonumber\\ \int_M \left(\int_M
\varphi^*(y) k_P^M(y,z,t)\,dy\right)f(z)\,dz=\int_M \varphi^*(z)
f(z)\,dz<\infty . \eea
 On the other hand, by the null-criticality, Fatou's lemma, and
(\ref{eqphiv}) we have
  \bean\label{eqcontra}
  \infty=\int_M \varphi^*(z) \alpha\,dz=
   \int_M \varphi^*(z) \lim_{n\to\infty}v(z,t_n)\,dz
\leq\\
 \liminf_{n\to\infty}\int_M \varphi^*(z) v(z,t_n)\,dz=\int_M
\varphi^*(z) f(z)\,dz<\infty. \eean Hence $\alpha=0$, and
therefore
\begin{equation}\label{eqlim0}
\lim_{t\to\infty}\int_M
k_P^M(x,y,t)f(y)\,dy=\lim_{t\to\infty}v(x,t)=0.
\end{equation}

Using the parabolic Harnack inequality and (\ref{eqinvar}), we
obtain that
\begin{equation}\label{eqHar}
k_P^M(x,y,t+t_n)\leq c_2(y)\varphi(x), \qquad {\and }\;
k_P^M(x,y,t+t_n)\leq c_1(x)\varphi^*(y)
\end{equation}
for all $x,y\in M$ and $t+t_n>1$ (see \cite{Pheat}).
 Now let $t_n\to
\infty$ be a sequence such that $\lim_{n\to\infty}
k_P^M(x,y,t+t_n)$ exists for all $(x,y,t)\in M\times M\times
\mathbb{R}$.  We denote the limit function by $u(x,y,t)$. It is
enough to show that any such $u$ is the zero solution. Recall that
as a function of $x$ and $t$, $u\in \mathcal{H}_+(M\times
\mathbb{R})$ (see \cite{Pheat}). Moreover, (\ref{eqHar}),  the
semigroup identity, and the dominated convergence theorem imply
that $$u(x,z,t+1)=\int_M u(x,y,t)k_P^M(y,z,1)\,dy.$$ It follows
that either $u=0$, or $u$ is a strictly positive function. On the
other hand, Fatou's lemma and (\ref{eqlim0}) imply that
 $$\int_M u(x,y,0)f(y)\,dy\leq \lim_{n\to
\infty}\int_M k_P^M(x,y,t_n)f(y)\,dy=0.$$
 Since  $f\gneqq 0$, it follows that $u=0$.\qed

Let $P$ be an elliptic operator of the form (\ref{P}) such that
$\lambda_0\geq 0$, and let $v\in \mathcal{C}_{P}(M)$ and $v^*\in
\mathcal{C}_{P^*}(M)$. It is well known \cite{PSt} that
\begin{equation}\label{eqinvar1}
\int_M k_P^M(x,y,t)v(y)\, dy\leq v(x),\;\;  \mbox{and } \int_M
k_P^M(x,y,t)v^*(x)\, dx\leq v^*(y).
\end{equation}
The parabolic Harnack inequality and (\ref{eqinvar1}) imply that
\begin{equation}\label{eqHar1}
k_P^M(x,y,t)\leq c_1(y)v(x), \qquad {\and }\; k_P^M(x,y,t)\leq
c_2(x)v^*(y)
\end{equation}
for all $x,y\in M$ and $t>1$ (see \cite{Pheat}). Recall that in
the critical case, $v$ and $v^*$ are in fact the ground states
$\varphi$ and $\varphi^*$ of $P$ and $P^*$ respectively, and by
(\ref{eqinvar}), we have equalities in (\ref{eqinvar1}).

 We now use theorems \ref{thm1} and \ref{mainthm}, estimate (\ref{eqHar1}),
and the dominated convergence theorem to strengthen Lemma
\ref{lemVa} for initial conditions which satisfy a certain
integrability condition.
 \begin{corollary}\label{Cormainthm} Let $P$ be an elliptic operator of
the form (\ref{P}) such that $\lambda_0\geq 0$. Let $f$ be a
continuous function on $M$ such that $v^*f\in L^1(M)$ for some
$v^*\in \mathcal{C}_{P^*}(M)$. Let
 $$u(x,t)=\int_{M} k_P^{M}(x,y,t)f(y)\,dy$$
 be the minimal solution of the Cauchy problem with initial data $f$ on
$M$. Fix $K\subset\subset M$.   Then $$\lim_{t\to\infty}\sup_{x\in
K}|u(x,t)-\mathcal{F}(x)|=0,$$ where $$\mathcal{F}(x)=
  \begin{cases}
    \varphi(x)\frac{\int_M \varphi^*(y)f(y)\,dy}{\int_M
\varphi^*(y)\varphi(y)\,dy} & \text{{\em if $P$ is
positive-critical in $M$}},
\\
    0 & \text{{\em otherwise}}.
  \end{cases}
$$
\end{corollary}

Suppose now that $P\mathbf{1}=0$ and $\int_M k_P^M(\cdot,y,t)\,
dy=\mathbf{1}$ (i.e. $\mathbf{1}$ is a positive invariant solution
of the operator $P$ in $M$). Corollary \ref{Cormainthm} implies
that for any $j\geq 1$ and all $x\in M$ we have
\begin{equation}\label{eqvMstar}
\lim_{t\to\infty}\int_{M_j^*} k_P^{M}(x,y,t)\,dy=
  \begin{cases}
    \frac{\int_{M_j^*} \varphi^*(y)\,dy}{\int_M
\varphi^*(y)\,dy} & \text{if $P$ is positive-critical in $M$},
\\
    1 & \text{otherwise}.
  \end{cases}\nonumber
\end{equation}
Therefore, if $P$ is not positive-critical in $M$, and $f$ is a
bounded continuous function such that $\liminf_{x\to\infty}
f(x)=\varepsilon>0$, then
\begin{equation}\label{eqliminf}
  \liminf_{t\to\infty}\int_{M} k_P^{M}(x,y,t)f(y)\,dy\geq \varepsilon.
\end{equation}
 Hence, if the integrability condition of Corollary \ref{Cormainthm}
is not satisfied, then the large time behavior of the minimal
solution of the Cauchy problem may be complicated. The following
example of W.~Kirsch and B.~Simon \cite{KS} demonstrates this
phenomenon.
\begin{example}\label{exKS}{\em
Consider the heat equation in $\mathbb{R}^d$. Let $R_j=e^{e^j}$
and let
 $$ f(x)=2+(-1)^j \qquad  \mbox{ if } \quad R_j<\sup_{1\leq i\leq d}
|y_i|<R_{j+1}, \;j\geq 1.$$
 Let $u$ be the minimal solution of the Cauchy problem with  initial data $f$.
Then for $t\sim R_jR_{j+1}$ one has that $u(0,t)\sim 2+(-1)^j$,
and thus $u(0,t)$ does not have a limit. Note that by Lemma
\ref{lemVa}, for $d=1$, $u(x,t)$ has exactly the same asymptotic
behavior as $u(0,t)$ for all $x\in \mathbb{R}$.
 }\end{example}

\mysection{Remarks and open problems}\label{secopen} In this
section, we mention some general open problems that are related to
the large time behavior of the heat kernel. The first conjecture
deals with the exact long time asymptotics of the heat kernel.

\begin{Conj}[E.~B.~Davies \cite{D}]\label{conjD}
Let $L=u_t+P(x, \partial_x)$ be a parabolic operator which is
defined on a Riemannian manifold $M$. Fix a reference point
$x_0\in M$. Then the limit
\begin{equation}\label{eqconjD}
\lim_{t\to\infty}\frac{k_P^M(x,y,t)}{k_P^M(x_0,x_0,t)}
\end{equation}
exists and is positive for all $x,y\in M$.
\end{Conj}
The answer to this conjecture seems to be closely related to the
question of the existence of a $\lambda_0$-invariant positive
solution (see \cite{D,PSt}).

The second conjecture was posed by the author \cite[Conjecture
3.6]{Pheat}.
\begin{Conj}\label{conjP}
Suppose that $P$ is a critical operator in $M$, then the ground
state $\varphi$ is a minimal positive solution in the cone
$\mathcal{H}_+(M\times \mathbb{R})$ of all nonnegative solutions
of the parabolic equation $Lu=0$ in $M\times \mathbb{R}$.
\end{Conj}
As noticed in \cite{Pheat}, if the conjecture is true, then
Theorem \ref{mainthm} would follow from (\ref{eqHar}).

 Recall also that by the parabolic Martin representation theorem,
the minimal positive solutions in $\mathcal{H}_+(M\times
\mathbb{R})$  are all parabolic Martin functions. Note that in the
positive-critical case, the ground state is clearly a parabolic
Martin function $K$ which corresponds to a fundamental sequence of
the form $\{(t_n,y_0)\}$, where $t_n\to -\infty$ and $y_0$ is a
fixed point in $M$. Indeed, by the definition of a Martin function
and Theorem \ref{thm1}, we have
\begin{equation}\label{eqMar}
K(x,y_0,t)=\lim_{n\to\infty}\frac{k_P^M(x,y_0,t-t_n)}{k_P^M(x_0,y_0,-t_n)}=
  \frac{\varphi(x)}{\varphi(x_0)}.
\end{equation}
On the other hand, if Conjecture \ref{eqconjD} is true, then it
can be easily checked that (\ref{eqMar}) is valid also in the
null-critical case and therefore, the ground state is always a
Martin function.

Recently K.~Burdzy and T.~S.~Salisbury \cite{BS} raised the
following more general problem
\begin{question}\label{questBS}
Determine which minimal harmonic functions are minimal in
$\mathcal{H}_+(M\times \mathbb{R})$, the cone of all parabolic
functions.
\end{question}
For more details see \cite{BS}.
\begin{center}
{\bf Acknowledgments} \end{center} The author wishes to thank
A.~Grigor'yan and R.~Pinsky for valuable discussions. This work
was partially supported by the Fund for the Promotion of Research
at the Technion.

  \end{document}